\documentclass{article}

\usepackage{amsmath}
\usepackage{amssymb}
\usepackage{amsthm}
\usepackage{braket}
\usepackage[colorlinks=true,urlcolor=blue,linkcolor=blue,plainpages=false,pdfpagelabels]{hyperref}
\usepackage{a4wide}

\newenvironment{claimproof}
{\begin{proof}[Proof of the claim]%
}
{\end{proof}}

\theoremstyle{plain}
\newtheorem{theorem}{Theorem}[section]

\newtheorem{lemma}[theorem]{Lemma}
\newtheorem{proposition}[theorem]{Proposition}
\newtheorem{remark}[theorem]{Remark}

\newcommand{\lp}{\left(}
\newcommand{\rp}{\right)}
\newcommand{\lb}{\left[}
\newcommand{\rb}{\right]}
\newcommand{\lbr}{\left\lbrace}
\newcommand{\rbr}{\right\rbrace}

\newcommand{\R}{{\mathbb R}}
\newcommand{\N}{{\mathbb N}}
\newcommand{\abs}[1]{\left\vert#1\right\vert}
\newcommand{\norm}[1]{\left\lVert#1\right\rVert}
\newcommand{\ceil}[1]{\left\lceil#1\right\rceil}
\newcommand{\Fix}{\operatorname{Fix}}

\newcommand{\rateoned}{\theta}
\newcommand{\ratetonemu}{\sigma}

\newcommand{\ratetwodelta}{\overline{\theta}}

\newcommand{\rateonealpha}{\lambda}
\newcommand{\rateonebeta}{\gamma_1}
\newcommand{\ratetwobeta}{\gamma_2}

\newcommand{\ratetwor}{\lambda_2}

\newcommand{\Mabd}{M}
\newcommand{\Mr}{M_r}
\newcommand{\zt}{z_t}
\newcommand{\zr}{z_r}
\newcommand{\zs}{z_s}
\newcommand{\zP}{z_P}
\newcommand{\zK}{z_K}

\newcommand{\kp}{K_p}
\newcommand{\kpz}{K_p^0}
\newcommand{\Kpasreg}{L_p}
\newcommand{\gammant}{\gamma_n^t}
\newcommand{\gammanr}{\gamma_n^r}
\newcommand{\gammans}{\gamma_n^s}

\newcommand{\gammanK}{\gamma_n^K}
\newcommand{\gammantr}{\gamma_n^{t_r}}
\newcommand{\gammanpt}{\gamma_{n+1}^t}
\newcommand{\gammanpr}{\gamma_{n+1}^r}
\newcommand{\gammanpP}{\gamma_{n+1}^P}
\newcommand{\gammanpK}{\gamma_{n+1}^K}
\newcommand{\psiB}{\psi_2}
\newcommand{\psiC}{\psi_1}
\newcommand{\ratemetaxn}{\Sigma}
\newcommand{\ratemetaxnr}{\Sigma^*}
\newcommand{\ynt}{y_{n,t}}

\title{Rates of metastability for a generalization of the Halpern iteration}
\author{Nicoleta Dumitru$^a$ and Lauren\c tiu Leu\c stean$^{a, b, c}$ \\ [2mm]
\footnotesize{$^a$ Faculty of Mathematics and Computer Science, University of Bucharest} \\
\footnotesize{$^b$ Simion Stoilow Institute of Mathematics of the Romanian Academy} \\
\footnotesize{$^c$ Institute for Logic and Data Science, Bucharest} \\
\footnotesize{E-mails: \href{mailto:nicoleta.dumitru@my.fmi.unibuc.ro}{nicoleta.dumitru@my.fmi.unibuc.ro}, \href{mailto:laurentiu.leustean@unibuc.ro}{laurentiu.leustean@unibuc.ro}}
}
\date{}

\begin{document}

\maketitle
	
\begin{abstract}
\noindent We apply proof mining methods developed by Kohlenbach and the second author to compute uniform rates of metastability for a generalization of the Halpern iteration studied by Kanzow and Shehu. As special cases, we obtain rates of metastability for the Halpern and modified Halpern iterations, as well as for the Tikhonov--Mann iteration. Furthermore,  convergence results are obtained  as immediate consequences of our quantitative results.
\end{abstract}

\medskip


\medskip


\section{Introduction}

Let $X$ be a Hilbert space, $T:X \to X$ a nonexpansive mapping with fixed points,  
$u\in X$, $(r_n)_{n\in \N}$ a sequence in $X$, and $(\alpha_n)_{n\in \N}$, $(\beta_n)_{n\in \N}$, 
$(\delta_n)_{n\in \N}$ be sequences in $[0,1]$ such that for all $n\in \N$,
\[\alpha_n + \beta_n + \delta_n \leq 1.\]

In this paper we study the asymptotic behavior of two iterations: the generalized Halpern iteration 
$(x_n)_{n\in \N}$, defined by 
\begin{equation}\label{def-gen-Halpern}
x_0 \in X, \qquad x_{n+1} = \delta_n u + \alpha_n x_n + \beta_n T x_n,
\end{equation}
and the inexact generalized Halpern iteration $(y_n)_{n\in\N}$, defined by
\begin{equation}\label{def-inexact-gen-Halpern}
y_0\in X,  \qquad y_{n+1}=\delta_nu+\alpha_ny_n+\beta_nTy_n+r_n.
\end{equation}
Kanzow and Shehu \cite[Section~4]{KanShe17} studied the inexact generalized Halpern iteration and proved a 
strong convergence result for this iteration. The well-known Halpern iteration is 
a special case of \eqref{def-gen-Halpern}: take $\alpha_n:=0$ and $\beta_n:=1-\delta_n$. 
Furthermore, the modified Halpern iteration, introduced by Kim and Xu \cite{KimXu05}, is another special case of \eqref{def-gen-Halpern}, as it is 
obtained from  \eqref{def-gen-Halpern} by letting 
$\alpha_n:= (1 - \delta_n)\alpha_n$ and $\beta_n:= (1 - \delta_n)(1 - \alpha_n)$ for all $n\in \N$.

\medskip 

The main quantitative results of the paper, proved in Section~\ref{sec-main-results},  provide 
uniform rates of metastability (in the sense of Tao \cite{Tao08a, Tao08}) for the generalized Halpern iteration and its inexact variant. These results are obtained 
by applying proof mining methods developed by Kohlenbach and the second author 
\cite{KohLeu12a,KohLeu12,KohLeu14} for the Halpern iteration. As special cases, we obtain 
rates of metastability for the Halpern and modified Halpern iterations, as well as for the Tikhonov--Mann iteration, defined by Cheval and the second 
author \cite{CheLeu22}. Furthermore, qualitative convergence results are obtained in 
Subsection~\ref{sec-main-results-strong-convergence} as immediate consequences of our 
quantitative results. 

Recall that a sequence $(z_n)_{n\in \N}$ in $X$  is Cauchy if and only if $(z_n)$ is metastable in 
the following sense: for all $k\in\N, g:\N\to\N$ there exists $N\in\N$ such that 
\[
\forall i,j\in\lb N;N+g(N)\rb
\,\, \lp \norm{z_i-z_j}\leq\frac1{k+1} \rp.
\]
A rate of metastability of $(z_n)$ is a mapping $\Omega:\N\times\N^\N\to\N$ satisfying
\[
\forall k\in\N\ \forall g:\N\to\N\ \exists N\leq\Omega(k,g) \forall i,j\in\lb N;N+g(N)\rb
\,\, \lp \norm{z_i-z_j}\leq\frac1{k+1} \rp.
\]
We refer to \cite{Koh19} for a discussion on (rates of) metastability. 

\medskip 

\textbf{Notations.} We denote $\N^* = \N \setminus \lbr 0 \rbr$.  If $m,n\in \N$, $m\leq n$, 
then we write $[m;n]=[m,n]\cap \N = \{m, m+1, \ldots, n\}$. 
Let $f : \N \to \N$. Then  $f^{(0)} (n) = n$ and $f^{(i+1)} (n) = f(f^{(i)} (n))$ for all $i, n\in\N$, and 
$\tilde{f}, f^+:\N \to \N$  are defined by 
\[\tilde{f}(n)=n+f(n), \quad f^+(n) = \max \lbr f(i) \mid 0 \leq i \leq n \rbr.\]
For $x\in \R$, $\ceil{x}$ denotes the ceiling of $x$. 
The set of fixed points of a mapping $T:X\to X$ is denoted by $\Fix(T)$.

\section{Quantitative lemmas} \label{sec-quant-res}

In what follows, we state several quantitative lemmas that will be used to prove our main quantitative 
theorems in Section~\ref{sec-main-results}. The first three lemmas are immediate. Lemma~\ref{quant-Xu-an-cn} 
is a quantitative version of a special case of a well-known lemma due to Xu \cite{Xu02}. Together with Lemma \ref{lemma-meta+conv-meta}, 
it will be used in the proof of Theorem~\ref{thm-main-metastability-inexact-genHalpern}. 
Lemma~\ref{lem52-KohLeu12-kp}, an adaptation of \cite[Lemma 5.2]{KohLeu12a} to our setting, is the main 
quantitative lemma used in the proof of Theorem~\ref{thm-main-metastability-genHalpern}. 

Let $(a_n)_{n \in \N}$ be a sequence of real numbers and $\varphi : \N \to \N$.  If $\lim\limits_{n\to \infty} a_n = a$, then 
$\varphi$ is a rate of convergence of $(a_n)$ to $a$ if 
\[\forall k \in \N \, \forall n \geq \varphi(k)  \lp \abs{a_n - a} \leq \frac{1}{k+1} \rp,\]
and  $\varphi$ is a Cauchy modulus of $(a_n)$ if 
\[\forall k \in \N \; \forall n \geq \varphi(k) \, \forall l \in \N  \lp \abs{a_{n + l} - a_n} \leq \frac{1}{k+1} \rp.\]
If $(a_n) \subseteq [0,\infty)$, $\varphi$ is a Cauchy modulus of 
$\sum\limits_{n=0}^\infty a_n$  if
\[\forall k\in \N\, \forall n\geq \varphi(k)\, \forall l\in\N^* \lp \sum\limits_{i=n+1}^{n+l} a_i\leq \frac1{k+1}\rp,\] 
and $\varphi$ is a  rate of divergence of $\sum\limits_{n=0}^\infty a_n$ if 
$\sum\limits_{i = 0}^{\varphi(n)} a_i \geq n$ for all $n\in \N$.

\begin{lemma}\cite[Lemma 2.3(i)]{FirLeu26a}\label{series-theta-n-1}\\
Let $(a_n)_{n \in \N} \subseteq [0,1]$ be such that $\sum\limits_{n=0}^\infty a_n$ diverges with rate of divergence $\theta$.
Then $\theta(n)\geq n-1$ for all $n\in\N$.
\end{lemma}

\begin{lemma}\label{series-Cauchy-rate}
Let $(a_n)_{n \in \N} \subseteq [0,\infty)$ be such that $\sum\limits_{n=0}^\infty a_n < \infty$ 
with Cauchy modulus $\varphi$. 
\begin{enumerate}
\item\label{series-Cauchy-rate-p1} $\forall k\in \N\, \forall n\geq \varphi(k)+1\, \forall l\in\N \lp \sum\limits_{i=n}^{n+l} 
a_i\leq \frac1{k+1}\rp$.
\item\label{series-Cauchy-ubd-M} $\sum\limits_{n=0}^\infty a_n \leq M$, where $M \in \N^*$ is such that 
$M \geq \sum\limits_{i = 0}^{\varphi(0)} a_i + 1$. 
\item\label{series-Cauchy-rate-p2} For all $c\in (0,\infty)$, $\sum\limits_{n=0}^\infty ca_n < \infty$ with 
Cauchy modulus $\psi(k)=\varphi(\ceil{c(k+1)}-1)$. 
\end{enumerate}
\end{lemma}

\begin{lemma}\label{basic-rate-conv-c}
Let $(a_n)_{n \in \N} \subseteq [0,\infty)$ be such that $\lim\limits_{n\to \infty}a_n =0 $  
with rate of convergence $\varphi$. Then 
for all $c\in (0,\infty)$,  $\lim\limits_{n\to \infty}(ca_n) =0 $ with rate of convergence 
$\psi(k)=\varphi\lp \lceil c(k+1)\rceil-1\rp$. 
\end{lemma}

\begin{lemma}\label{quant-Xu-an-cn} \cite[Proposition 2.7]{LeuPin24} \\
Let $(s_n)_{n \in \N}, (c_n)_{n \in \N}\subseteq [0,+\infty)$  and $(a_n)_{n \in \N}\subseteq [0,1]$ satisfy, for all $n\in\N$,
\[s_{n+1}\leq (1-a_n)s_n + c_n. \]
Assume that $L\in\N^*$ is an upper bound on $(s_n)$, $\sum\limits_{n=0}^\infty a_n$  diverges with rate of 
divergence $\theta$, and  $\sum\limits_{n=0}^\infty c_n$ converges 
with Cauchy modulus $\chi$.

Then $\lim\limits_{n\to\infty} s_n=0$ with rate of convergence $\Sigma$ defined by 
\[
\Sigma(k)=\theta \lp \chi(2k+1)+1+ \ceil{\ln(2L(k+1))}\rp +1.
\]
\end{lemma}

\begin{lemma}\label{lemma-meta+conv-meta}\cite[Proposition 7]{LeuPin21}\\
Let $(u_n)_{n\in \N}, (v_n)_{n\in \N}$ be sequences in $X$ such that $(u_n)$ is 
Cauchy with rate of metastability $\Omega$ and $\lim\limits_{n\to\infty} \norm{u_n-v_n}=0$ with rate of convergence $\varphi$.

Then $(v_n)$ is Cauchy with rate of metastability $\Gamma$ given by 
\[
\Gamma(k,g)=\max\{\varphi(3k+2),\ \Omega(3k+2,\widehat{g}_k)\},
\]
with $\widehat{g}_k(n)=\max\{\varphi(3k+2), n\}-n+g(\max\{\varphi(3k+2), n\})$.
\end{lemma}

\begin{lemma}\label{lem52-KohLeu12-kp}
Let $(\lambda_n)_{n\in \N}$, $(b_n)_{n\in \N}$, $(c_n)_{n\in \N}$, $(s_n)_{n\in \N}$ be 
real sequences, $M \in \N^*$, and mappings $\theta, \psiC:\N \to \N$ with the following properties:
\begin{enumerate}
\item[(H1)] $(s_n)$ is bounded with upper bound $M$; 
\item[(H2)] $(\lambda_n)\subseteq [0,1]$ and $\sum\limits_{n = 0}^{\infty} \lambda_n = \infty$ with 
rate of divergence $\theta$; 
\item[(H3)] $(c_n)\subseteq [0,\infty)$ and $\sum\limits_{n=0}^\infty c_n<\infty$ with Cauchy modulus $\psiC$.
\end{enumerate}
Assume, furthermore, that $l\in \N$, $g : \N \to \N$, and $\psiB:\N \to \N$ are such that 
\begin{enumerate}
\item[(H4)] for all $n \geq \psiB(l)$, $b_n \leq \frac1{6(l+1)}$;
\item[(H5)] for all $n \in \N$,
\begin{equation}\label{ineq-main-quant-lemma}
s_{n+1} \leq (1 - \lambda_n) s_n + \lambda_n b_n + c_n + \frac1{\Gamma(l,g)},
\end{equation}
where 
\begin{align*}
\Gamma:\N \times \N^\N \to \N, & \quad \Gamma(n,h) = 3(n+1)(\Theta(n) + h(\Theta(n)) - \psi(n)),\\
\Theta:\N\to \N, & \quad \Theta(n) =  \theta \lp \psi(n) + \ceil{\ln (3M (n+1))} \rp + 1, \\
\psi:\N\to \N, & \quad \psi(n) = \max\lbr \psiB(n), \psiC(6n+5)+1\rbr.
\end{align*}
\end{enumerate}
Then
$$
 \forall n \in [\Theta(l);\Theta(l) + g(\Theta(l))] \; \lp s_n \leq \frac{1}{l+1} \rp.
$$
\end{lemma}
\begin{proof}
Let $l \in \N$ and $g : \N \to \N$. Set
\begin{equation}\label{def-H-main-quant}
 \quad H_l = \Theta(l) - \psi(l) - 1 = \theta\lp \psi(l) + \ceil{\ln \lp 3M(l+1) \rp} \rp - \psi(l).
\end{equation}
As $(\lambda_n)\subseteq [0,1]$ and $\theta$ is a rate of divergence of 
$\sum\limits_{n = 0}^{\infty} \lambda_n$, by Lemma~\ref{series-theta-n-1}, we have
\begin{align*}
H_l & \geq \ceil{\ln \lp 3M(l+1) \rp} - 1 \geq \ceil{\ln 3} - 1 \geq 1.
\end{align*}
As $\Gamma(l,g)=3(l+1)(H_l+1 + g(\Theta(l)))$, it is easy to see that $\Gamma(l,g)\geq 6$.

As in the proof of \cite[Lemma 5.2]{KohLeu12a}, we obtain the following inequalities: for all $m \in \N$ and all $n \geq \psiB(l)$,
\begin{align}
s_{n+m+1} & \leq M \prod\limits_{j = n}^{n+m} (1 - \lambda_j) + \lb 1 - \prod\limits_{j = n}^{n+m}(1-\lambda_j) \rb \frac{1}{6(l+1)} + 
\sum\limits_{j = n}^{n+m}c_j + \frac{m+1}{\Gamma(l,g)}, \label{lem-snmp-1} \\
s_{n+m+1}  & \leq M \exp\lp -\sum\limits_{j = n}^{n+m} \lambda_j \rp + \frac{1}{6(l+1)} +  
\sum\limits_{j = n}^{n+m}c_j + \frac{m+1}{\Gamma(l,g)}. \label{lem-snmp-2} 
\end{align}

\newcounter{quantno}

\textbf{Claim \stepcounter{quantno}\arabic{quantno}:} For all $m \in \N$ and all $n \geq \psi(l)$,
\[
s_{n+m+1} \leq M \exp\lp -\sum\limits_{j = n}^{n+m} \lambda_j \rp + \frac{1}{3(l+1)} +  \frac{m+1}{\Gamma(l,g)}.
\]
\begin{claimproof} Since $n \geq \psi(l)\geq \psiB(l)$, we can apply \eqref{lem-snmp-2} to get that
\[
s_{n+m+1} \leq M \exp\lp -\sum\limits_{j = n}^{n+m} \lambda_j \rp + \frac{1}{6(l+1)} +  \sum\limits_{j = n}^{n+m}c_j + \frac{m+1}{\Gamma(l,g)}.
\]
As $n \geq \psi(l)\geq \psiC(6l+5)+1$ and, by hypothesis, $\psiC$ is a Cauchy modulus of $\sum\limits_{n=0}^\infty c_n$, 
we can apply Lemma \ref{series-Cauchy-rate}\eqref{series-Cauchy-rate-p1} to get that $\sum\limits_{j = n}^{n+m}c_j \leq \frac{1}{6(l+1)}$. 
The conclusion follows. 
\end{claimproof}

\textbf{Claim \stepcounter{quantno}\arabic{quantno}:} 
For all $m \geq H_l$, 
\[
M \exp\lp - \sum\limits_{j=\psi(l)}^{\psi(l)+m} \lambda_j  \rp \leq \frac{1}{3(l+1)}.
\]
\begin{claimproof} 
As in the proof of \cite[Lemma 5.2]{KohLeu12a} (see also the proof of \cite[Lemma 9]{Leu07a}), we get that for all
$m\in \N$, $n\in \N^*$, 
\begin{equation}\label{lem-snmp-3-eq11}
M \exp\lp - \sum\limits_{j=n}^{n+m} \lambda_j \rp \leq \frac{1}{3(l+1)}  
 \quad \text{iff}\quad \sum\limits_{j=0}^{n+m} \lambda_j \geq \sum\limits_{j=0}^{n-1} \lambda_j + \ln \lp 3M(l+1) \rp.
\end{equation}
As for all $m \geq H_l$,
\begin{align*}
\sum\limits_{j=0}^{\psi(l)+m} \lambda_j &\geq \sum\limits_{j=0}^{\psi(l)+H_l} \lambda_j 
\stackrel{\eqref{def-H-main-quant}}{=} \sum\limits_{j=0}^{\Theta(l)-1} \lambda_j  
\geq \psi(l) + \ceil{\ln (3M (l+1))} 
 \geq \sum\limits_{j = 0}^{\psi(l) - 1} \lambda_j + \ceil{\ln(3M(l+1))},
\end{align*}
we can apply \eqref{lem-snmp-3-eq11} with $n := \psi(l)\in \N^*$ to get the conclusion. 
\end{claimproof}

\textbf{Claim \stepcounter{quantno}\arabic{quantno}:} 
For all $m \geq H_l$, 
\[
s_{\psi(l)+m+1} \leq \frac{2}{3(l+1)} + \frac{m+1}{\Gamma(l,g)}.
\]
\begin{claimproof} 
Apply Claim 1 with $n:=\psi(l)$ and Claim 2 to get that for all $m \geq H_l$,
\begin{align*}
s_{\psi(l)+m+1} & \leq M \exp\lp -\sum\limits_{j = \psi(l)}^{\psi(l) +m} \lambda_j \rp + \frac{1}{3(l+1)} 
+ \frac{m+1}{\Gamma(l,g)} \leq \frac{2}{3(l+1)} + \frac{m+1}{\Gamma(l,g)}.
\end{align*}
\end{claimproof}

We now conclude the proof. Let $n \in [\Theta(l);\Theta(l) + g(\Theta(l))]$. Then 
\begin{align*}
s_n & = s_{\psi(l)+(n-\psi(l)-1)+1} \leq  \frac{2}{3(l+1)} + \frac{n - \psi(l)}{\Gamma(l,g)} \\
& \quad \text{ by Claim 3 with $m:=n-\psi(l)-1 \geq \Theta(l) - \psi(l) -1 = H_l$}\\
& \leq \frac{2}{3(l+1)} + \frac{\Theta(l) - \psi(l) + 
g(\Theta(l))}{\Gamma(l,g)} \quad \text{as $n \leq \Theta(l) + g(\Theta(l))$} \\
&  \leq \frac1{l+1} \quad \text{by the definition of $\Gamma$}.
\end{align*}

\end{proof}

\section{A Browder-type sequence}

Let $X$ be a real Hilbert space, $u\in X$, and $T:X\to X$ be a nonexpansive mapping with $\Fix(T)\ne \emptyset$. For any $t\in (0,1)$, define $S_t:X\to X$ by 
$S_t(x)=tu+(1-t)Tx$. It is well known that $S_t$ is a $(1-t)$-contraction, 
hence one can apply Banach Contraction Mapping Principle to get a unique fixed point $\zt\in X$. Thus, 
for all $t\in (0,1)$,
\begin{equation}\label{def-zt}
\zt=tu+(1-t)T\zt,
\end{equation}
and, as a consequence, 
\begin{equation}\label{zt-u-minus-Tzt}
\zt - T \zt  = t (u - T\zt), \quad u - \zt  = (1-t)(u - T\zt).
\end{equation}

\begin{lemma} \label{ubd-lemma-main-zt}
Let $p \in \Fix(T)$. Then  for all $t\in (0,1)$,
\begin{enumerate}
\item\label{ubd-p-zt-Tzt} $\norm{p - T\zt} \leq \norm{p - \zt} \leq \norm{u - p}$;
\item\label{ubd-zt} $\norm{\zt} \leq \norm{u - p} + \norm{p}$;
\item\label{ubd-uzt-Tzt} $\norm{u - \zt}, \norm{u - T \zt}  \leq 2 \norm{u - p}$.
\end{enumerate}
\end{lemma} 
\begin{proof}
\begin{enumerate}
\item Since $T$ is nonexpansive, $\norm{p - T\zt}= \norm{Tp - T\zt}\leq \norm{p - \zt}$. Furthermore, 
\begin{align*}
\norm{p - \zt} & \stackrel{\eqref{def-zt}}{=} \norm{p - (t u + (1 - t) T\zt)} 
 = \norm{(1 - t) (p - T\zt) + t (p - u)} \\
& \leq (1 - t) \norm{p - T\zt} + t \norm{p - u} \leq (1 - t) \norm{p - \zt} + t \norm{u - p}. 
\end{align*}
Thus,  $t\norm{p - \zt}  \leq t \norm{u - p}$, so $\norm{p - \zt} \leq \norm{p - u}$. 
\item $\norm{\zt}  \leq  \norm{\zt-p} + \norm{p} \stackrel{\eqref{ubd-p-zt-Tzt}}{\leq} \norm{u - p} + \norm{p}$.
\item We have  $\norm{u - \zt} \leq \norm{u - p} + \norm{p - \zt} 
\stackrel{\eqref{ubd-p-zt-Tzt}}{\leq}   2\norm{u - p}$, and $\norm{u - T\zt} \leq \norm{u - p} + \norm{p - T\zt} 
\stackrel{\eqref{ubd-p-zt-Tzt}}{\leq}   2\norm{u - p}$.
\end{enumerate}
\end{proof}

The following proposition, which plays a key role in establishing our main results, follows from the quantitative version of Browder's theorem proved by
Kohlenbach~\cite[Theorem~4.2]{Koh11}.

\begin{proposition} \label{Browder-quant-genHalpern}
Assume that $(t_n)_{n\in \N}$ is a nonincreasing sequence in $(0,1)$. Then $(z_{t_n})$ is Cauchy with rate of metastability 
$\Omega$ given by
\begin{equation*}
\Omega(k,g)=\tilde{g}^{\lp 4M^2(k+1)^2 \rp}(0),
\end{equation*}
where $M\in \N^*$ is such that $M\geq \norm{u - p}$ for some $p\in \Fix(T)$. 
\end{proposition}
\begin{proof}
By Lemma \ref{ubd-lemma-main-zt}\eqref{ubd-uzt-Tzt}, 
\[\norm{z_{t_n}-u} \leq 2M\]
 for all $n\in \N$. 
Apply \cite[Theorem~4.2]{Koh11} with  $U:=T$,  $v_0:=u$, $s_n:=1-t_n$, $\check{u}_n:=z_{t_n}$, $d:=2M$, and 
$\varepsilon:=\frac1{k+1}$.
\end{proof}

\section{Generalized Halpern iteration}

Let $X$ be a real Hilbert space, $T:X \to X$ a nonexpansive mapping with fixed points, $u\in X$, and 
$(\alpha_n)_{n\in \N}$, $(\beta_n)_{n\in \N}$, $(\delta_n)_{n\in \N}$ be sequences in $[0,1]$ such that for all $n\in \N$,
\[
\alpha_n + \beta_n + \delta_n \leq 1.\]

For all $n\in\N$, denote $\mu_n= 1-\alpha_n - \beta_n - \delta_n$. Obviously, 
$\mu_n \in [0,1]$ for all $n\in \N$.

The generalized Halpern iteration $(x_n)_{n\in \N}$ is defined by \eqref{def-gen-Halpern}:
\[ x_0 \in X, \qquad x_{n+1} = \delta_n u + \alpha_n x_n + \beta_n T x_n.\]

Let us recall the following well-known inequalities that will be used in the sequel:

\begin{lemma}
\begin{enumerate}
\item For all $x, y \in X$, 
\begin{equation}
\norm{x + y}^2 \leq \norm{x}^2 + 2 \braket{y, x+y}. \label{xpy2-leq-x2pbryxpy}
\end{equation}
\item For all $x,y,z\in X$ and $a,b,c\in[0,1]$ satisfying $a+b+c\leq 1$, 
\begin{equation}
\|ax+by+cz\|^2\leq a\|x\|^2+b\|y\|^2+c\|z\|^2.\label{eq-Hilbert-useful-abc}
\end{equation}
\end{enumerate}
\end{lemma}

\begin{lemma} \label{lem-zt-xn-1}
For all $n \in \N$ and $t\in (0,1)$,
\begin{align}
x_n - \zt & = (1 - t) (x_n - T\zt) + t(x_n - u), \label{xn-minus-zt-1} \\
\norm{x_n - \zt}^2 & \leq (1 - t)^2 \norm{x_n - T\zt}^2 + 2t\norm{x_n - \zt}^2 - 2t\braket{u - \zt, x_n - \zt}, \label{norm-xn-zt-sq-1}  \\
\norm{T x_n - \zt}^2 & \leq \norm{x_n - \zt}^2 + t^2 \norm{u - T \zt}^2 + 2 t \norm{x_n - \zt} 
\norm{u - T \zt}, \label{norm-Txn-minus-zt-2} \\
2 t \braket{u - \zt, x_n - \zt} & \leq (1 - t)^2 \norm{x_n - Tx_n} \lp \norm{x_n - Tx_n} + 2\norm{x_n - \zt}\rp + t^2\norm{x_n - \zt}^2. 
\label{useful-ineq-gammant-1}
\end{align}
\end{lemma}
\begin{proof}
 \eqref{xn-minus-zt-1} is an 
immediate application of \eqref{def-zt}: 
$$x_n - \zt = x_n - (t u + (1 - t) T \zt)= (1 - t) (x_n - T\zt) + t (x_n - u).$$

The inequalities \eqref{norm-xn-zt-sq-1} and \eqref{norm-Txn-minus-zt-2} are obtained as follows:
\begin{align*}
\norm{x_n - \zt}^2 & \stackrel{\eqref{xn-minus-zt-1}}{=} 
\norm{\lp 1 - t \rp (x_n - T\zt) + t (x_n - u)}^2 \\
& \stackrel{\eqref{xpy2-leq-x2pbryxpy}}{\leq}  
\norm{(1 - t)(x_n - T\zt)}^2 + 2 \braket{t (x_n - u), x_n - \zt} \\
& = (1 - t)^2 \norm{x_n - T\zt}^2 + 2 t \braket{x_n - u, x_n - \zt}\\
& = (1 - t)^2 \norm{x_n - T\zt}^2 + 2t\lp \braket{(x_n - \zt) + (\zt - u), x_n - \zt}\rp \\
& = (1 - t)^2 \norm{x_n - T\zt}^2 + 2t\norm{x_n - \zt}^2 - 2t\braket{u - \zt, x_n - \zt} \\[1mm]
\norm{Tx_n - \zt}^2 &  = \norm{T x_n - T \zt}^2 + \norm{T \zt - \zt}^2 + 2 \braket{T x_n - T \zt, T \zt - \zt} \\
& \leq \norm{T x_n - T \zt}^2 + \norm{T \zt - \zt}^2 + 2 \norm{T x_n - T \zt} \norm{T \zt - \zt}  \\
& \leq \norm{x_n - \zt}^2 + \norm{T \zt - \zt}^2 + 2 \norm{x_n - \zt} \norm{T \zt - \zt} \quad \text{since $T$ is nonexpansive} \\
& \stackrel{\eqref{zt-u-minus-Tzt}}{=} \norm{x_n - \zt}^2 + t^2 \norm{u - T \zt}^2 + 2 t \norm{x_n - \zt} \norm{u - T \zt}.
\end{align*}
Note that 
\begin{align*}
\norm{x_n - \zt}^2 & \stackrel{\eqref{norm-xn-zt-sq-1}}{\leq}  (1 - t)^2 \norm{x_n - T\zt}^2 + 2 t \norm{x_n - \zt}^2 - 2 t \braket{u - \zt, x_n - \zt} \\
& = (1 - t)^2 \lp \norm{x_n - Tx_n}^2 + \norm{Tx_n - T\zt}^2 + 2 \braket{x_n - Tx_n,Tx_n - T\zt} \rp + 2 t \norm{x_n - \zt}^2 \\
& \quad - 2 t \braket{u - \zt, x_n - \zt} \\
& \leq (1 - t)^2 \norm{x_n - Tx_n}^2 + (1 - t)^2 \norm{Tx_n - T\zt}^2 + 2 (1 - t)^2
 \norm{x_n - Tx_n} \norm{Tx_n - T\zt} \\
& \quad + 2 t \norm{x_n - \zt}^2 - 2 t \braket{u - \zt, x_n - \zt} \\
& \leq (1 - t)^2 \norm{x_n - Tx_n}^2 + (1 - t)^2 \norm{x_n - \zt}^2 + 2 (1 - t)^2
 \norm{x_n - Tx_n} \norm{x_n - \zt} \\
& \quad + 2 t \norm{x_n - \zt}^2 - 2 t \braket{u - \zt, x_n - \zt} \qquad  \text{since $T$ is nonexpansive}\\
& = (1 - t)^2 \norm{x_n - Tx_n} \lp \norm{x_n - Tx_n} + 2\norm{x_n - \zt}\rp + (1+t^2) \norm{x_n - \zt}^2 \\ 
& \quad - 2 t \braket{u - \zt, x_n - \zt}, 
 \end{align*}
hence \eqref{useful-ineq-gammant-1} follows.
\end{proof}

Denote, for all $n \in \N$ and $t\in (0,1)$,
\begin{equation}\label{def-ynt}
\ynt  = \norm{u - T \zt}^2   + 2 \norm{x_n - \zt} \norm{u - T \zt}.
\end{equation} 
 
\begin{lemma}
For all $n \in \N$ and $t\in (0,1)$,
\begin{equation}\label{norm-xnzt2-ubd-ab}
\alpha_n \norm{x_n - \zt}^2 +  \beta_n \norm{T x_n - \zt}^2  \leq 
(1-\delta_n)\norm{x_n - \zt}^2 +(1-\alpha_n-\delta_n)t \ynt.
\end{equation} 
\end{lemma} 
\begin{proof}
Let us denote, for simplicity, $a_n= \alpha_n \norm{x_n - \zt}^2 +  \beta_n \norm{T x_n - \zt}^2 $. We get that
\begin{align*}
a_n & \stackrel{\eqref{norm-Txn-minus-zt-2}}{\leq} \alpha_n \norm{x_n - \zt}^2 +  \beta_n \lp \norm{x_n - \zt}^2 + 
t^2 \norm{u - T \zt}^2 + 2 t \norm{x_n - \zt} \norm{u - T \zt} \rp\\
& = (\alpha_n+\beta_n)\norm{x_n - \zt}^2 + \beta_n  \lp t^2\norm{u - T \zt}^2 + 
2 t \norm{x_n - \zt} \norm{u - T \zt} \rp \\
& \leq  (1-\delta_n)\norm{x_n - \zt}^2 + (1-\alpha_n-\delta_n)t \lp t\norm{u - T \zt}^2 
+ 2  \norm{x_n - \zt} \norm{u - T \zt}\rp \\
& \quad \text{as $\alpha_n + \beta_n\leq 1- \delta_n$ and $\beta_n \leq 1-\alpha_n-\delta_n$}\\
& \leq  (1-\delta_n)\norm{x_n - \zt}^2 + (1-\alpha_n-\delta_n)t\ynt, \quad \text{as $t \in (0,1)$}.
\end{align*}
\end{proof}

\begin{lemma}
For all $n \in \N$ and $t\in (0,1)$, 
\begin{align}
x_{n+1} - \zt & = \alpha_n (x_n - \zt) + \beta_n (T x_n - \zt) + \delta_n(u - \zt) - \mu_n\zt, 
\label{xnp1-minus-zt-1}\\
\norm{x_{n+1} - \zt}^2 & \leq  (1-\delta_n)\norm{x_n - \zt}^2 
+2\delta_n \braket{u - \zt, x_{n+1} - \zt} + \mu_n \norm{\zt}^2 + (1-\alpha_n-\delta_n)t\ynt. \label{norm-xnp-minus-zt2-2}
\end{align}
\end{lemma}
\begin{proof}
\eqref{xnp1-minus-zt-1} follows directly from the definition of $(x_n)$. Furthermore,
\begin{align*}
\norm{x_{n+1} - \zt}^2 & \stackrel{\eqref{xnp1-minus-zt-1}}{=}  
\norm{\alpha_n (x_n - \zt) + \beta_n (T x_n - \zt) + \delta_n(u - \zt) - \mu_n\zt}^2 \\
& \leq \norm{\alpha_n (x_n - \zt) +  \beta_n (T x_n - \zt) - \mu_n\zt}^2 + 
2 \braket{\delta_n(u - \zt), x_{n+1} - \zt} \\
& \quad \text{by \eqref{xpy2-leq-x2pbryxpy} with } x:=\alpha_n (x_n - \zt) + \beta_n (T x_n - \zt)  - \mu_n\zt
\text{ and } y:=\delta_n(u - \zt)  \\
& \leq \alpha_n\norm{x_n - \zt}^2 + \beta_n\norm{T x_n - \zt}^2 
+ \mu_n \norm{\zt}^2+ 2\delta_n \braket{u - \zt, x_{n+1} - \zt} \\
& \quad \text{by \eqref{eq-Hilbert-useful-abc}, as $\alpha_n+\beta_n+\mu_n=1-\delta_n \leq 1$.}
\end{align*}
Now apply \eqref{norm-xnzt2-ubd-ab} to get \eqref{norm-xnp-minus-zt2-2}.

\end{proof}

For the rest of the section, we assume that $\sum\limits_{n=0}^\infty \mu_n<\infty$ with an upper bound
$\Mabd\in \N$.

Let $p\in \Fix(T)$ and take $\kpz, \kp\in \N^*$ as follows: 
\begin{align} 
\kpz  \geq \max \lbr \norm{x_0 - p}, \norm{u - p}, \norm{p} \rbr, \qquad
\kp  = (1 + \Mabd) \kpz. \label{def-kzp-kp} 
\end{align}

\begin{lemma}  \label{ubd-lemma-xnzt-main}
For all $n \in \N$ and $t\in (0,1)$,
\begin{enumerate}
\item \label{ubd-xn-p} $\norm{Tx_n - p} \leq \norm{x_n - p} \leq \kp$;
\item \label{ubd-xn-Txn-zt} $\norm{x_n - Tx_n}, \norm{x_n - \zt}, \norm{u - T \zt} \leq 2\kp$.
\end{enumerate}
\end{lemma} 
\begin{proof}
\begin{enumerate}
\item The fact that $\norm{x_n - p} \leq \kp$ is an immediate consequence of \cite[Lemma~3.5(ii)]{DumLeu26}. 
Obviously, $\norm{Tx_n - p} =  \norm{Tx_n - Tp}\leq \norm{x_n - p}$. 
\item We have  $\norm{x_n-Tx_n} \leq  \norm{x_n - p} + \norm{p - Tx_n} 
\leq 2\kp$. Furthermore, by Lemma \ref{ubd-lemma-main-zt}\eqref{ubd-p-zt-Tzt},\eqref{ubd-uzt-Tzt}
we obtain
\begin{align*}
\norm{x_n - \zt} & \leq \norm{x_n - p} + \norm{p - \zt} 
 \leq \norm{x_n - p} + \norm{p - u}  \leq  \kp + \kpz \leq 2 \kp, \\
 \norm{u - T \zt} & \leq 2 \norm{u-p}\leq 2\kpz\leq 2\kp.
\end{align*}
\end{enumerate}
\end{proof}

As an immediate consequence of Lemma \ref{ubd-lemma-xnzt-main}\eqref{ubd-xn-Txn-zt} and of the definition \eqref{def-ynt} of $\ynt$, we get 
that for all $n \in \N$ and $t\in (0,1)$,
\begin{equation}\label{ubd-ytn-12kp2}
\ynt \leq 12\kp^2.
\end{equation}

\subsection{An auxiliary sequence}

For all $n \in \N$ and $t\in (0,1)$, define 
\begin{equation} \label{def-gammant}
\gammant= 2 \braket{u - \zt, x_n - \zt} - 4t \kp^2.
\end{equation}

\begin{lemma}
For all $n \in \N$ and $t\in (0,1)$,
\begin{equation}\label{xnzt-main-ineq}
\norm{x_{n+1} - \zt}^2  \leq 
(1 - \delta_n) \norm{x_n - \zt}^2 + \delta_n \gammanpt + \mu_n \norm{\zt}^2  + 12t \kp^2.
\end{equation}
\end{lemma}
\begin{proof}
We have
\begin{align*}
\norm{x_{n+1} - \zt}^2 & \stackrel{\eqref{norm-xnp-minus-zt2-2}}{\leq}   (1-\delta_n)\norm{x_n - \zt}^2 
+2\delta_n \braket{u - \zt, x_{n+1} - \zt} + \mu_n \norm{\zt}^2 + (1-\alpha_n-\delta_n)t\ynt \\[2mm]
& \leq  (1-\delta_n)\norm{x_n - \zt}^2 
+ 2\delta_n \braket{u - \zt, x_{n+1} - \zt} + \mu_n \norm{\zt}^2  + (1-\delta_n)t\ynt \\
& \quad \text{as $1-\alpha_n-\delta_n \leq 1-\delta_n$}\\
& \stackrel{\eqref{ubd-ytn-12kp2}}{\leq} (1-\delta_n)\norm{x_n - \zt}^2
+ 2\delta_n \braket{u - \zt, x_{n+1} - \zt}  + \mu_n \norm{\zt}^2  + 12(1-\delta_n)t\kp^2 \\
& \stackrel{\eqref{def-gammant}}{=}  (1-\delta_n)\norm{x_n - \zt}^2 + \delta_n \gammanpt +  \mu_n \norm{\zt}^2 +
4\delta_n t \kp^2  + 12(1-\delta_n)t\kp^2\\
& = (1 - \delta_n) \norm{x_n - \zt}^2 + \delta_n \gammanpt + \mu_n \norm{\zt}^2  + 4t \kp^2 + 8(1 - \delta_n)t\kp^2\\
& \leq (1 - \delta_n) \norm{x_n - \zt}^2 + \delta_n \gammanpt+ \mu_n \norm{\zt}^2  + 12t \kp^2 \quad 
\text{since $\delta_n \geq 0$}.
\end{align*}
\end{proof}

\begin{lemma}
For all $n \in \N$ and $t\in (0,1)$,
\begin{equation}
\gammant  \leq \frac{6\kp(1 - t)^2}{t} \norm{x_n - Tx_n}.\label{gammant-leq-xnutxn}
\end{equation}
\end{lemma}
\begin{proof}
We have 
\begin{align*}
2t\braket{u - \zt, x_n - \zt} & \stackrel{\eqref{useful-ineq-gammant-1}}{\leq}   (1 - t)^2 \norm{x_n - Tx_n} \lp \norm{x_n - Tx_n} + 
2\norm{x_n - \zt}\rp + t^2\norm{x_n - \zt}^2 \\
& \leq   6\kp(1 - t)^2 \norm{x_n - Tx_n} + 4\kp^2t^2 \quad 
\text{by Lemma \ref{ubd-lemma-xnzt-main}\eqref{ubd-xn-Txn-zt}}.
\end{align*}
The conclusion follows. 
\end{proof}

\section{Main results} \label{sec-main-results}

Throughout this section, $X$ is a real Hilbert space, $T:X \to X$ is a nonexpansive mapping with $\Fix(T)\neq\emptyset$, 
$u\in X$, and $(\alpha_n)_{n\in \N}$, $(\beta_n)_{n\in \N}$, 
$(\delta_n)_{n\in \N}$,  $(\mu_n)_{n\in \N}$ are sequences in $[0,1]$ such that for all $n\in \N$,
\[\alpha_n + \beta_n + \delta_n \leq 1, \qquad \mu_n= 1-\alpha_n - \beta_n - \delta_n.\]
Let $(x_n)_{n\in\N}$ be the generalized Halpern iteration defined by \eqref{def-gen-Halpern}:
\[
x_0 \in X,  \qquad x_{n+1} = \delta_n u + \alpha_n x_n + \beta_n T x_n. 
\]

The main quantitative result of the paper provides a uniform rate of metastability for $(x_n)$, assuming that a rate of $T$-asymptotic regularity of $(x_n)$ is available. 
Recall that $(x_n)$ is said to be $T$-asymptotically regular if 
$\lim\limits_{n \to \infty} \norm{x_n - T x_n} = 0$. A rate of convergence of the sequence 
$(\norm{x_n - T x_n})$ to $0$ will be 
referred to as a rate of $T$-asymptotic regularity of $(x_n)$. 

\begin{theorem}\label{thm-main-metastability-genHalpern}
Assume that the following hold:
\begin{enumerate}
\item\label{mtmg-Tasreg} $(x_n)$ is $T$-asymptotically regular with rate $\Psi$; 
\item\label{mtmg-mun} $\sum\limits_{n=0}^\infty \mu_n < \infty$ with Cauchy modulus 
$\ratetonemu$ and an upper bound $\Mabd \in \N$;
\item\label{mtmg-deltan} $\sum\limits_{n=0}^\infty \delta_n = \infty$ with rate of divergence $\rateoned$;
\item\label{mtmg-Kp} $\kp$ is defined by \eqref{def-kzp-kp}  for some $p\in \Fix(T)$.
\end{enumerate}
Then $(x_n)$ is Cauchy with rate of metastability $\ratemetaxn$, that is, for all $k \in \N$ and all $g : \N \to \N$, 
$$
\exists N \leq \ratemetaxn(k, g) \; \forall m, n \in [N;N + g(N)] \; \lp \norm{x_n - x_m} \leq \frac{1}{k+1} \rp,
$$
with 
\begin{equation}\label{def-ratemetaxn}
\ratemetaxn: \N \times \N^\N \to \N, \quad 
\ratemetaxn(k,g) =  \rateoned^+ \lp \Lambda_{k,g} + \ceil{\ln\lp 48\kp^2(k+1)^2\rp}\rp + 1, 
\end{equation}
where for all $r\in \N$, $\Omega_r, \psi_r, \Theta_r:\N \to \N$ and $\Gamma_r:\N \times \N^\N \to \N$  are defined as follows:
\begin{align*}
\Omega_r(n) & = \Psi(6\kp(r+1)(n+1)-1), \\
\psi_r(n) & = \max\lbr \Omega_r(12n+11)-1, \ratetonemu(24\kp^2(n+1)-1)+1\rbr,\\
\Theta_r(n) & = \rateoned\lp \psi_r(n) + \ceil{\ln (12\kp^2 (n+1))} \rp + 1, \\
\Gamma_r(n,h) & = 3(n+1) (\Theta_r(n) +h(\Theta_r(n)) - \psi_r(n)),
\end{align*}
and for all $k\in \N$ and all $g:\N\to \N$, $k_0,\Pi_{k,g},\Pi^*_{k,g}\in \N$ and 
$f_{k,g},f_{k,g}^*:\N\to \N$ are defined as follows:
\begin{align*}
k_0 & = 192\kp(\kp+2)(k+1)^2-1, \\
f_{k,g}(n)  & = \max \lbr 12\kp^2\Gamma_n(4(k+1)^2-1,g)-n-1, 0 \rbr,  \\
f_{k,g}^*(n) & = f_{k,g}(n+k_0) + k_0,\\
\Pi_{k,g}  &= \widetilde{f_{k,g}^*}^{\lp 4\kp^2(k_0+1)^2 \rp}(0), \\
\Pi^*_{k,g}  & = \Pi_{k,g} + k_0, \\
\Lambda_{k,g} & = \max\lbr \psi_l(4(k+1)^2-1)  \mid  k_0 \leq l\leq \Pi^*_{k,g} \rbr.
\end{align*}
\end{theorem}
\begin{proof}
Let $k \in \N$ and $g : \N \to \N$ be fixed. Consider, for all $r\in \N$, $t_r=\frac1{r+1}$, and denote, for simplicity, 
$z_{t_r}$  by $z_r$ and $\gammantr$ by $\gammanr$. 

\medskip

\newcounter{stepno}

\textbf{Claim \stepcounter{stepno}\arabic{stepno}} For all $r\in \N$ and  all $n \geq \Omega_r(k)$, 
$\gammanr \leq \frac1{k+1}$.
\begin{claimproof} Let $r\in \N$ and $n\geq \Omega_r(k)=\Psi(6\kp(r+1)(k+1)-1)$. As 
$\Psi$ is a rate of $T$-asymptotic regularity of $(x_n)$, we have 
\begin{align*}
\norm{x_n - Tx_n} & \leq  \frac1{6\kp(r+1)(k+1)}.
\end{align*} 
Now apply \eqref{gammant-leq-xnutxn} to obtain
\begin{align*}
\gammanr & \leq \frac{6\kp(1 - t_r)^2}{t_r} \norm{x_n - Tx_n} \leq \frac{6\kp}{t_r} \norm{x_n - Tx_n} 
= 6\kp(r+1)\norm{x_n - Tx_n} \leq \frac1{k+1}. 
\end{align*} 
\end{claimproof}

\textbf{Claim \stepcounter{stepno}\arabic{stepno}}  For all $r,n\in \N$, 
\begin{equation}\label{xnzt-main-ineq-main-r}
\norm{x_{n+1} - \zr}^2  \leq 
(1 - \delta_n) \norm{x_n - \zr}^2 + \delta_n \gammanpr + 4\kp^2\mu_n   + \frac{12\kp^2}{r+1}.
\end{equation}
\begin{claimproof} 
We have
\begin{align*}
\norm{x_{n+1} - \zr}^2  & \leq 
(1 - \delta_n) \norm{x_n - \zr}^2 + \delta_n \gammanpr + \mu_n \norm{\zr}^2  + \frac{12\kp^2}{r+1} \quad \text{by \eqref{xnzt-main-ineq} with $t:= t_r=\frac1{r+1}$}\\
& \leq (1 - \delta_n) \norm{x_n - \zr}^2 + \delta_n \gammanpr + 4\kp^2\mu_n   + \frac{12\kp^2}{r+1},
\end{align*}
as, by Lemma~\ref{ubd-lemma-main-zt}\eqref{ubd-zt} and the definition \eqref{def-kzp-kp} of $\kpz$, $\kp$, 
$\norm{\zr} \leq \norm{u - p} + \norm{p} \leq 2\kpz \leq  2\kp$. 
\end{claimproof}

\textbf{Claim \stepcounter{stepno}\arabic{stepno}} For all  $r,s,n\in\N$, 
\begin{equation}
\gammanr - \gammans  \leq 8\kp\norm{\zs - \zr} + \frac{4\kp^2}{s + 1} - \frac{4\kp^2}{r + 1}.
\label{def-gammanr-gammans}
\end{equation}
\begin{claimproof}
We have
\begin{align*}
\gammanr - \gammans &  \stackrel{\eqref{def-gammant}}{=} 2 \braket{u - \zr, x_n - \zr} - 4 t_r \kp^2 - 
\lp 2 \braket{u - \zs, x_n - \zs}  - 4 t_s \kp^2 \rp \\
& = 2 \braket{u - \zr, x_n - \zr} - 2 \braket{u - \zs, x_n - \zs} + 4 t_s\kp^2 -  4 t_r \kp^2 \\
& = 2 \braket{u - \zr, x_n - \zr} - 2\braket{(u - \zr) + (\zr - \zs), x_n - \zs}  + \frac{4\kp^2}{s + 1} - \frac{4\kp^2}{r + 1}\\
& = 2 \braket{u - \zr, \zs - \zr} - 2\braket{\zr - \zs, x_n - \zs} + \frac{4\kp^2}{s + 1} - \frac{4\kp^2}{r + 1} \\
& \leq 2 \norm{\zs - \zr} \lp \norm{u - \zr} + \norm{x_n - \zs} \rp  + \frac{4\kp^2}{s + 1} - \frac{4\kp^2}{r + 1} \\
& \leq 8\kp\norm{\zs - \zr} + \frac{4\kp^2}{s + 1} - \frac{4\kp^2}{r + 1} 
\quad \text{by Lemmas~\ref{ubd-lemma-main-zt}\eqref{ubd-uzt-Tzt} and \ref{ubd-lemma-xnzt-main}\eqref{ubd-xn-Txn-zt}.}
\end{align*}
\end{claimproof}

\textbf{Claim \stepcounter{stepno}\arabic{stepno}} There exists $K_1\in \N$, $K_1 \leq \Pi_{k,g}$ such that 
\begin{equation}\label{ineq-useful-K1}
\forall q,s \in [K_1;K_1+f_{k,g}^*(K_1)]  \lp \norm{z_q-z_s} \leq \frac{1}{k_0+1} \rp.
\end{equation}
\begin{claimproof} As $(t_n)$ is nonincreasing,  we can apply Proposition~\ref{Browder-quant-genHalpern} 
with $M:=\kp \geq \norm{u - p}$ to get that 
\[
\forall l \in \N \, \forall h:\N \to \N \, \exists K_1 \leq \Omega(l,h) \, \forall q,s \in 
[K_1;K_1+h(K_1)] \lp \norm{z_q-z_s} \leq \frac{1}{l+1} \rp,
\]
where $\Omega(l,h)=\tilde{h}^{(4\kp^2(l+1)^2)}(0)$. 
In particular, this holds for $l:=k_0$ and $h:=f_{k,g}^*$. Then 
\begin{align*}
\Omega(k_0,f_{k,g}^*) & =  \widetilde{f_{k,g}^*}^{\lp 4\kp^2(k_0+1)^2 \rp}(0) 
= \Pi_{k,g}.
\end{align*}
\end{claimproof}

\textbf{Claim \stepcounter{stepno}\arabic{stepno}} Let $K = K_1 + k_0$. Then 
\begin{equation}\label{ineq-useful-K0}
K\in [k_0;\Pi^*_{k,g}]  \quad \text{and} \quad \forall q,s \in [K; K+f_{k,g}(K)]  \lp \norm{z_q-z_s} \leq \frac{1}{k_0+1} \rp.
\end{equation}
\begin{claimproof} Since, by the previous claim, $K_1 \leq \Pi_{k,g}$, we get that 
$k_0 \leq K \leq \Pi_{k,g} + k_0 = \Pi^*_{k,g}$. Furthermore, $K_1 \leq K$ and 
$K + f_{k,g}(K) = K_1 + k_0 + f_{k,g}(K_1 + k_0) = K_1 + f_{k,g}^*(K_1)$, hence 
$[K;K + f_{k,g}(K)] \subseteq [K_1;K_1+f_{k,g}^*(K_1)]$.
Now apply \eqref{ineq-useful-K1}.  
\end{claimproof}

\textbf{Claim \stepcounter{stepno}\arabic{stepno}} For all  $r\in[K;K+f_{k,g}(K)]$ and all $n\in\N$, 
\[ \gammanr - \gammanK  \leq \frac1{48(k+1)^2}.\]
\begin{claimproof} Let  $r\in[K;K+f_{k,g}(K)]$. We have, for all $n\in\N$,
\begin{align*}
\gammanr - \gammanK  & \stackrel{\eqref{def-gammanr-gammans}}{\leq} 8\kp\norm{\zK - \zr} + \frac{4\kp^2}{K + 1} - 
\frac{4\kp^2}{r + 1} 
\leq 8\kp\norm{\zK - \zr} +\frac{4\kp^2}{K + 1} \\
& \leq \frac{8\kp}{k_0+1} +  \frac{4\kp^2}{K + 1}  \quad \text{by Claim 5 applied with $q:=K$, $s:=r$} \\
& \leq \frac{4\kp(\kp+2)}{k_0+1} \quad \text{as $K\geq k_0$}\\
& = \frac1{48(k+1)^2} \quad \text{by the definition of $k_0$.} 
\end{align*}
\end{claimproof}

\textbf{Claim \stepcounter{stepno}\arabic{stepno}} For all $r\in[K;K+f_{k,g}(K)]$ and all $n\geq \Omega_K(48(k+1)^2-1)-1$, 
\[\gammanpr  \leq \frac1{24(k+1)^2}.\]
\begin{claimproof}
As $n+1\geq \Omega_K(48(k+1)^2-1)$, we can apply Claim 1 with $r:= K$, $n:=n+1$, and $k:=48(k+1)^2-1$ to get that 
$\gammanpK \leq \frac1{48(k+1)^2}$. By Claim 6, we conclude that 
$\gammanpr \leq  \gammanpK + \frac1{48(k+1)^2} \leq \frac1{24(k+1)^2}$.	
\end{claimproof}

\textbf{Claim \stepcounter{stepno}\arabic{stepno}} Let $P=K+f_{k,g}(K)$.  Then for all $n\in \N$, 
\begin{equation}\label{xnzt-main-ineq-main}
\norm{x_{n+1} - \zP}^2  \leq 
(1 - \delta_n) \norm{x_n - \zP}^2 + \delta_n \gammanpP + 4\kp^2\mu_n  + \frac1{\Gamma_K(4(k+1)^2-1,g)}.
\end{equation}
\begin{claimproof}
First note that $\Gamma_r(m,h)\in \N^*$ for all $r,m\in \N$ and all $h:\N\to \N$. 
 Apply Claim 2 with $r:=P$ to obtain that for all $n\in \N$, 
\begin{align*}
\norm{x_{n+1} - \zP}^2  & \leq 
(1 - \delta_n) \norm{x_n - \zP}^2 + \delta_n \gammanpP + 4\kp^2\mu_n   + \frac{12\kp^2}{P+1}.
\end{align*}
As 
\begin{align*}
P +1 & = f_{k,g}(K) + K + 1 = \max \lbr 12\kp^2\Gamma_K(4(k+1)^2-1,g), K +1\rbr \\
& \geq 12\kp^2\Gamma_K(4(k+1)^2-1,g),
\end{align*}
we get that 
\[
\frac{12\kp^2}{P+1} \leq  
\frac1{\Gamma_K(4(k+1)^2-1,g)}.\]
Then \eqref{xnzt-main-ineq-main} follows. 
\end{claimproof}

Let us set
\begin{equation}\label{def-N}
N=\Theta_K(4(k+1)^2-1).
\end{equation}

\textbf{Claim \stepcounter{stepno}\arabic{stepno}} The following holds: 
\begin{equation}\label{ineq-final-1-Theta-K}
\forall n \in [N;N+g(N)] \; \lp \norm{x_n - \zP}  \leq \frac{1}{2(k+1)} \rp.
\end{equation}
\begin{claimproof}
We apply Lemma~\ref{lem52-KohLeu12-kp} with the following choices:
\begin{align*}
s_n:=\norm{x_n - \zP}^2, \quad  \lambda_n:=\delta_n, \quad b_n:=\gammanpP, \quad c_n:=4\kp^2\mu_n, 
\quad M:=4\kp^2, \quad \theta:= \rateoned, 
\end{align*}
and 
\begin{align*}
 \psiC(n) := \ratetonemu(4\kp^2(n+1)-1), \quad \psiB(n):= \max\lbr\Omega_K(12n+11)-1,0\rbr, \quad l:=4(k+1)^2-1.
\end{align*}
Let us verify the hypotheses (H1)--(H5) of Lemma~\ref{lem52-KohLeu12-kp}.

(H1) holds by Lemma~\ref{ubd-lemma-xnzt-main}\eqref{ubd-xn-Txn-zt}, since $\norm{x_n - \zP}^2\leq 4\kp^2$ for all $n\in \N$.  
Hypothesis (H2) follows directly 
from the assumptions on $(\delta_n)$. To verify (H3), use the assumption that 
$\sum\limits_{n=0}^\infty \mu_n  < \infty$ with Cauchy modulus $\ratetonemu$
and apply Lemma~\ref{series-Cauchy-rate}\eqref{series-Cauchy-rate-p2}. Claim 7 gives  
$\gammanpP \leq \frac1{24(k+1)^2}=\frac1{6(l+1)}$ for all 
$n \geq \Omega_K(48(k+1)^2-1)-1$. Since $\psiB(l)=\psiB(4(k+1)^2-1)\geq \Omega_K(48(k+1)^2-1)-1$, we get (H4). 
Finally, in the notation of Lemma~\ref{lem52-KohLeu12-kp}, we have $\psi=\psi_K$, 
$\Theta=\Theta_K$, and $\Gamma=\Gamma_K$. Then (H5) follows from \eqref{xnzt-main-ineq-main}.

We may therefore apply Lemma~\ref{lem52-KohLeu12-kp} to obtain
\eqref{ineq-final-1-Theta-K}. 
\end{claimproof}

Let us conclude the proof. By \eqref{ineq-final-1-Theta-K} and the triangle inequality,
\[
\forall m,n\in [N;N+g(N)] \; \lp \norm{x_m-x_n} \leq \frac{1}{k+1} \rp.
\]
Furthermore, 
\begin{align*}
N & = \Theta_K(4(k+1)^2-1) = \rateoned \lp  \psi_K(4(k+1)^2-1)  + \ceil{\ln\lp 48\kp^2(k+1)^2\rp} \rp + 1 \\
& \leq \rateoned^+\lp \psi_K(4(k+1)^2-1)  + \ceil{\ln\lp 48\kp^2(k+1)^2\rp} \rp + 1.
\end{align*}
As $k_0 \leq K\leq \Pi^*_{k,g}$, we get that $\psi_K(4(k+1)^2-1)  \leq \Lambda_{k,g}$. 
Finally, as $\rateoned^+$ is increasing, we get that 
\[N\leq \rateoned^+ \lp \Lambda_{k,g} + \ceil{\ln \lp 48\kp^2(k+1)^2\rp}\rp + 1 = \ratemetaxn(k,g).\]

Thus, there exists $N\leq\ratemetaxn(k,g)$ such that
\[
\forall m, n \in [N;N + g(N)] \; \lp \norm{x_n - x_m} \leq \frac{1}{k+1} \rp.
\]
\end{proof}

\begin{remark}
In hypothesis \eqref{mtmg-mun} of Theorem~\ref{thm-main-metastability-genHalpern}, it suffices to assume
that the series $\sum\limits_{n=0}^{\infty}\mu_n$ converges with Cauchy
modulus $\ratetonemu$, as, by Lemma~\ref{series-Cauchy-rate}\eqref{series-Cauchy-ubd-M}, one may take
$$\Mabd = \ceil{\sum\limits_{i=0}^{\ratetonemu(0)} \mu_i} + 1.$$
We nevertheless formulate the hypothesis with $\Mabd$ as a separate parameter, since this 
allows us to take $\Mabd=0$ rather than $\Mabd=1$ when $\mu_n=0$ for all $n\in \N$.
\end{remark}

\medskip

The first hypothesis of Theorem~\ref{thm-main-metastability-genHalpern} is that $(x_n)$ is 
$T$-asymptotically regular with rate $\Psi$. Such rates were recently computed by the authors in \cite{DumLeu26}.
The following theorem is a consequence of Theorem~\ref{thm-main-metastability-genHalpern} together with   
the main result of \cite{DumLeu26}.

\begin{theorem} \label{thm-main-metastability-asreg-genHalpern}
Assume that the following hold: 
\begin{enumerate}
\item $\sum\limits_{n=0}^\infty \mu_n < \infty$ with Cauchy modulus $\ratetonemu$ and an upper 
bound $\Mabd \in \N$;
\item $\sum\limits_{n=0}^\infty \delta_n = \infty$ with rate of divergence $\rateoned$ and  
$\sum\limits_{n=0}^{\infty} \abs{\delta_n-\delta_{n+1}} <\infty$ with Cauchy modulus $\ratetwodelta$;
\item $\sum\limits_{n=0}^{\infty} \abs{\alpha_n-\alpha_{n+1}} <\infty$ with Cauchy modulus 
$\rateonealpha$;
\item $\sum\limits_{n=0}^{\infty} \abs{\beta_n-\beta_{n+1}} <\infty$ with Cauchy modulus $\rateonebeta$ 
and $\lim\limits_{n\to\infty} \beta_n=1$ with rate of convergence $\ratetwobeta$;
\item for some $p\in \Fix(T)$, $\kpz$, $\kp$ are defined by \eqref{def-kzp-kp}, and 
$\Kpasreg = (2 + \Mabd) \kpz + 1$. 
\end{enumerate} 
Then $(x_n)$ is Cauchy with rate of metastability $\ratemetaxn$ given by \eqref{def-ratemetaxn}, 
with $\Psi$ defined as follows:
$$\Psi(k)=\max \lbr \Phi(3k+2), \ratetwobeta\lp 6\Kpasreg(k+1) - 1 \rp\rbr,$$
where 
\begin{align*}
\Phi(k) & = \rateoned^+\lp \chi(2k+1)+ \ceil{\ln(4\Kpasreg(k+1))} + 2\rp, \\
\chi(k) & = \max\lbr\ratetwodelta(6\Kpasreg(k+1)-1), \rateonealpha(6\Kpasreg(k+1)-1),
\rateonebeta(6\Kpasreg(k+1)-1)\rbr.
\end{align*}
\end{theorem}
\begin{proof}
Applying \cite[Theorem~4.1]{DumLeu26} with $\rho:=0$ and $r_n:=0$ for all $n\in \N$ 
(see also \cite[Section~4.1]{DumLeu26}) yields a rate of $T$-asymptotic regularity. 
The result now follows from Theorem \ref{thm-main-metastability-genHalpern}.
\end{proof}

If we take $\alpha_n:=0$ and $\beta_n:=1-\delta_n$, then $\mu_n=0$ and the generalized Halpern iteration 
reduces to the classical Halpern iteration
\[ x_{n+1}=\delta_n u+(1-\delta_n)Tx_n.\]
Thus, rates of metastability for the Halpern iteration can be obtained as immediate consequences 
of our theorems. In Theorem~\ref{thm-main-metastability-genHalpern} we take  
$\Mabd=0$, $\kp=\kpz$, and $\ratetonemu \equiv 0$. Furthermore, in Theorem~\ref{thm-main-metastability-asreg-genHalpern}, 
we take
 $\Kpasreg = 2\kpz + 1$, 
$\rateonealpha \equiv 0$, $\rateonebeta = \ratetwodelta$, and $\ratetwobeta$ to be a rate of convergence of $(\delta_n)$ to $0$.
Hence, $\chi(k) = \ratetwodelta(6\Kpasreg(k+1)-1)$. 

By taking in the generalized Halpern iteration  \eqref{def-gen-Halpern} 
$\alpha_n:= (1 - \delta_n)\alpha_n$ and $\beta_n:= (1 - \delta_n)(1 - \alpha_n)$, we obtain the modified Halpern iteration, introduced by Kim and Xu \cite{KimXu05}:
\[x_{n + 1} = \delta_n u + (1 - \delta_n)(\alpha_n x_n + (1 - \alpha_n)T x_n).\]
Hence, we get rates of metastability for this iteration, too, by applying our main results. 
Such rates for the modified Halpern iteration in CAT(0) spaces were obtained by Schade and 
Kohlenbach \cite{SchKoh12} using similar methods. As pointed out by Cheval, Kohlenbach, and the 
second author \cite{CHeKOhLeu23}, from the rates of metastability for the modified Halpern iteration 
we get such rates for the Tikhonov--Mann iteration, defined by Cheval and the second 
author \cite{CheLeu22}, and, in particular, for the modified Mann iteration, studied 
in \cite{YaoZhoLio09, BotCseMei19}.

\medskip

The next theorem computes a rate of metastability for the inexact generalized Halpern iteration defined 
by \eqref{def-inexact-gen-Halpern}:
\[
y_0\in X,  \qquad y_{n+1}=\delta_nu+\alpha_ny_n+\beta_nTy_n+r_n,
\]
where $(r_n)_{n\in \N}$ is a sequence in $X$, assuming that a rate of metastability is available 
for the generalized Halpern iteration $(x_n)$ with $x_0=y_0$.

\begin{theorem}\label{thm-main-metastability-inexact-genHalpern}
Assume that the following hold:
\begin{enumerate}
\item $\ratemetaxn$ is a rate of metastability for the generalized Halpern iteration with $x_0=y_0$;
\item $\sum\limits_{n=0}^\infty \delta_n = \infty$ with rate of divergence $\rateoned$; 
\item $\sum\limits_{n=0}^{\infty} \norm{r_n} <\infty$ with Cauchy modulus $\ratetwor$ 
and an upper bound $\Mr\in \N^*$. 
\end{enumerate}
Then $(y_n)$ is Cauchy with rate of metastability $\ratemetaxnr$ defined by 
$$\ratemetaxnr: \N \times \N^\N \to \N, \quad  \ratemetaxnr(k,g) = \max\lbr \varphi(3k+2), \ratemetaxn(3k+2,\widehat{g}_k)\rbr,$$
where 
\begin{equation}\label{def-rate-vp-inexact-Halpern}
\varphi:\N\to \N, \quad \varphi(n)=\rateoned \lp \ratetwor(2n+1)+1+ \ceil{\ln(2\Mr(n+1))}\rp +1,
\end{equation}
and, for all $k\in \N$, $g:\N\to \N$, 
$$\widehat{g}_k:\N\to \N, \quad\widehat{g}_k(n) = \max\lbr \varphi(3k+2), n\rbr-n+g(\max\lbr \varphi(3k+2), n\rbr).$$
\end{theorem}
\begin{proof}
First, note that for all $n\in \N$, 
\begin{align*}
 \norm{x_{n+1}-y_{n+1}} & =  \norm{\alpha_n(x_n-y_n)+\beta_n(Tx_n-Ty_n) - r_n} 
 \le (\alpha_n+\beta_n)\norm{x_n-y_n}+\|r_n\| \\
 & \le (1-\delta_n)\norm{x_n-y_n}+\|r_n\|,
\end{align*}
so, in particular, $\norm{x_{n+1}-y_{n+1}} \leq \norm{x_n-y_n}+\|r_n\|$.  Hence, we get by induction on $n$ 
that $\norm{x_n-y_n} \leq \sum\limits_{i=0}^{n-1} \norm{r_i}\leq \Mr$ for all $n\in \N$.

We can apply Lemma \ref{quant-Xu-an-cn} with 
\[s_n:= \norm{x_n-y_n}, \quad a_n:=\delta_n,  \quad c_n:=\norm{r_n},\quad L:=\Mr, \quad 
\theta:=\rateoned, \quad\chi:=\ratetwor \]
to get that  $\lim\limits_{n\to\infty} \norm{x_n-y_n}=0$ with rate of convergence $\varphi$ defined by 
\eqref{def-rate-vp-inexact-Halpern}.

Finally, an application of Lemma \ref{lemma-meta+conv-meta} with $u_n:= x_n$, $v_n:=y_n$, 
$\Omega:=\ratemetaxn$ yields that $(y_n)$ is Cauchy with rate of metastability $\ratemetaxnr$.
\end{proof}

Theorem~\ref{thm-main-metastability-inexact-genHalpern} yields a similar result
for inexact versions of the Halpern and modified Halpern iterations.

\subsection{Convergence results} \label{sec-main-results-strong-convergence}

As is typically the case with quantitative results obtained using proof mining methods, disregarding 
the quantitative information yields new qualitative consequences. We now use the results obtained 
above to derive convergence results for the generalized Halpern iteration $(x_n)$ and 
the inexact generalized Halpern iteration $(y_n)$.

\begin{theorem}\label{thm-main-conv-genHalpern-asreg}
Assume that $(x_n)$ is $T$-asymptotically regular, $\sum\limits_{n=0}^\infty \mu_n < \infty$, and 
$\sum\limits_{n=0}^\infty \delta_n = \infty$. Then $(x_n)$ converges to a fixed point of $T$.
\end{theorem} 
\begin{proof}
By Theorem~\ref{thm-main-metastability-genHalpern}, $(x_n)$ is Cauchy, hence, 
since $X$ is a Hilbert space, $(x_n)$ converges to some $p\in X$. Since $(x_n)$ is 
$T$-asymptotically regular, it follows that $p$ is a fixed point of $T$. 
\end{proof}
 
\begin{theorem} \label{thm-main-conv-genHalpern}
Assume that the following hold:
\begin{enumerate}
\item $\sum\limits_{n=0}^\infty \delta_n=\infty$;
\item the series $\sum\limits_{n=0}^\infty\mu_n$, \, 
$\sum\limits_{n=0}^{\infty} \abs{\delta_n-\delta_{n+1}}$, \, 
$\sum\limits_{n=0}^{\infty} \abs{\alpha_n-\alpha_{n+1}}$, \, 
$\sum\limits_{n=0}^{\infty} \abs{\beta_n-\beta_{n+1}}$ converge;
\item $\lim\limits_{n\to\infty} \beta_n=1$.
\end{enumerate}
Then $(x_n)$ converges to a fixed point of $T$.
\end{theorem} 
\begin{proof}
By the proof of Theorem~\ref{thm-main-metastability-asreg-genHalpern}, $(x_n)$ is
$T$-asymptotically regular, so Theorem~\ref{thm-main-conv-genHalpern-asreg} applies.
\end{proof}

\begin{theorem} \label{thm-main-conv-inexact-genHalpern}
Under the hypotheses of Theorem~\ref{thm-main-conv-genHalpern}, assume, moreover, that 
$\sum\limits_{n=0}^{\infty} \norm{r_n}< \infty$.

Then $(y_n)$ converges and $\lim\limits_{n\to\infty} y_n = \lim\limits_{n\to\infty} x_n$. 
\end{theorem} 
\begin{proof}
The conclusion follows from Theorem~\ref{thm-main-conv-genHalpern} and the fact,
established in the proof of Theorem~\ref{thm-main-metastability-inexact-genHalpern}, 
that $\lim\limits_{n\to\infty} \norm{x_n-y_n}=0$.
\end{proof}

By taking $\alpha_n:=0$ and $\beta_n:=1-\delta_n$ in Theorem~\ref{thm-main-conv-inexact-genHalpern}, we get 
a convergence result for the following  
inexact version of the Halpern iteration:
\[y_{n+1}=\delta_n u+(1-\delta_n)Ty_n+r_n. \]

Kanzow and Shehu also proved, using arguments due to Maing\'e~\cite{Mai08}, the convergence 
of  the inexact generalized Halpern iteration $(y_n)$ (see \cite[Theorem~4.1]{KanShe17}) under the following assumptions:
\[
\lim\limits_{n\to\infty} \delta_n=0, \qquad \sum\limits_{n=0}^\infty \delta_n = \infty, \qquad
\liminf\limits_{n\to\infty} \alpha_n\beta_n > 0, \qquad \sum\limits_{n=0}^\infty\mu_n < \infty, \qquad 
\sum\limits_{n=0}^{\infty} \norm{r_n} < \infty.
\]
However, the assumption $\liminf\limits_{n\to\infty} \alpha_n\beta_n > 0$ excludes the Halpern 
iteration as a special case.  \\

\mbox{}

\textbf{Acknowledgements:} The second author used ChatGPT 5.6 (OpenAI) during the preparation of this paper for proofreading, 
including the correction of typographical and grammatical errors and improvements to English phrasing. The authors take full responsibility for the content of the paper.

\end{document}